\documentclass[11pt,reqno]{amsart}
\usepackage[T1]{fontenc}
\usepackage{lmodern}
\usepackage{amsmath,amssymb,mathtools}
\usepackage{microtype}
\usepackage[letterpaper,margin=1.08in,headheight=14pt]{geometry}
\usepackage{enumitem}
\usepackage{booktabs}
\usepackage{xcolor}
\usepackage[colorlinks=true,linkcolor=blue!40!black,citecolor=blue!40!black,urlcolor=blue!45!black]{hyperref}
\usepackage{bookmark}
\hypersetup{pdftitle={Uniform ultrafilters of maximal character at aleph omega},pdfsubject={Expanded reconstruction of Gitik's argument with an ordinary-ultrafilter extension},pdfauthor={Reconstructed research draft}}
\newtheorem{theorem}{Theorem}

\newtheorem{lemma}[theorem]{Lemma}

\theoremstyle{definition}

\theoremstyle{remark}

\numberwithin{equation}{section}
\DeclareMathOperator{\cf}{cf}

\newcommand{\uU}{\mathfrak u}

\newcommand{\Bcal}{\mathcal B}

\setlist[enumerate]{label=\textup{(\roman*)},leftmargin=2.2em,itemsep=3pt}
\setlist[itemize]{leftmargin=1.8em,itemsep=3pt}
\title[Large ultrafilter number]{A large ultrafilter number at $\aleph_\omega$}
\author{Fanxin Wu}
\address{Department of Mathematics, University of California, Irvine, Irvine, CA 92697}

\date{\today}
\email{fanxin.wu.math@gmail.com}
\date{\today}

\begin{document}
\begin{abstract}
Gitik constructed a model where $\aleph_\omega$ is strong limit, $2^{\aleph_\omega}=\aleph_{\omega+2}$ and every strongly uniform ultrafilter on $\aleph_\omega$ has character $\aleph_{\omega+2}$. We show that in fact every uniform ultrafilter on $\aleph_\omega$ has character $\aleph_{\omega+2}$ in that model.
\end{abstract}
\maketitle

\section{Introduction}\label{sec:intro}

An ultrafilter $D$ on an infinite set $X$ is \emph{uniform} if $|A|=|X|$ for every $A\in D$. For an ideal $I$ on $X$ and $A,B\subseteq X$, write $A\subseteq_I B$ if $A\setminus B\in I$. A \emph{base of $D$ modulo $I$} is a family $\Bcal\subseteq D$ such that for every $A\in D$ there is some $B\in\Bcal$ such that $B\subseteq_I A$. Set
\[
 \chi_I(D)=\min\{|\Bcal|:\Bcal\text{ is a base of }D\text{ modulo }I\},\qquad
 \uU_I(\kappa)=\min\{\chi_I(D):D\text{ is uniform on }\kappa\}.
\]
So when $I$ is the trivial ideal $\{\emptyset\}$ these recover the usual character $\chi(D)$ and ultrafilter number $\uU(\kappa)$, and in general we have $\chi_I(D)\leq\chi(D)$ and $\uU_I(\kappa)\leq\uU(\kappa)$.

We shall mainly be concerned with the case when $X$ is a singular cardinal $\kappa$ and $I=[\kappa]^{<\kappa}$; here we also denote $\chi_I(D),\uU_I(\kappa)$ by $\chi^*(D),\uU^*(\kappa)$ respectively. At one point we shall also consider the slightly more general case $I=[X]^{<\kappa}$ for some set $X$; here we also denote $\chi_I(D)$ by $\chi_{<\kappa}(D)$.

Garti and Shelah \cite{garti-shelah2012ultrafilter} initiated the study of ultrafilter numbers at singular cardinal; in particular they showed the consistency of $\kappa^+=\uU^*(\kappa)=\uU(\kappa)<2^\kappa$ relative to a supercompact cardinal. Garti, Gitik and Shelah \cite{garti-gitik-shelah2020cardinal} achieved this at $\aleph_\omega$ and lowered the large cardinal hypothesis to the optimal one, a measurable cardinal $\kappa$ with $o(\kappa)=\kappa^{++}$. There have since been many further results on $\uU(\kappa)$ for singular $\kappa$ \cite{benhamou-jirattikansakul2023small,honzik-stejskalova2022small,usuba2026monotonicity}.

A natural question is whether $\uU^*(\kappa)>\kappa^+$ or even $\uU(\kappa)>\kappa^+$ is consistent for a strong limit singular $\kappa$. Gitik \cite{gitik2019spectrum} made important progress in this direction: starting from a model of $\mathsf{GCH}$ with a $(\kappa,\kappa^{++})$-extender, he constructed a model that satisfies the following:
\begin{enumerate}
\item $\aleph_\omega$ is a strong limit and $2^{\aleph_\omega}=\aleph_{\omega+2}$;
\item every strongly uniform ultrafilter $D$ on $\aleph_\omega$ satisfies $\chi^*(D)=\aleph_{\omega+2}$.
\end{enumerate}
Here an ultrafilter $D$ on a singular cardinal $\kappa$ is called \emph{strongly uniform} if there is a strictly increasing cofinal sequence of regular cardinals $\langle\lambda_i:i<\cf(\kappa)\rangle$ below $\kappa$ such that
\[
 \{i<\cf(\kappa):|A\cap\lambda_i|=\lambda_i\} \text{ is unbounded in }\cf(\kappa) \qquad(A\in D).
\]

The exact value of $\uU^*(\aleph_\omega)$ in this model was left open, and Gitik suggested that it might be $\aleph_{\omega+1}$. Here we shall show that in fact it is $\aleph_{\omega+2}$, thus answering the above question in the positive. This is a simple consequence of two facts: a dichotomy about uniform ultrafilters on $\aleph_\omega$ and a lower bound for the character of a regular ultrafilter under suitable pcf hypotheses.

\begin{lemma}\label{lem:dichotomy}
Let $D$ be a uniform ultrafilter on $\kappa=\aleph_\omega$. Then either $D$ is $(\omega,\kappa)$-regular, or $D$ is Rudin--Keisler above a strongly uniform ultrafilter on $\kappa$.
\end{lemma}

Recall that $D$ is $(\theta,\kappa)$-regular if there is a family of $\kappa$ many sets in $D$ such that the intersection of any $\theta$ many of them is empty.

\begin{lemma}\label{lem:regular-bound}
Suppose $\langle\lambda_i:i<\omega\rangle$ are strictly increasing regular cardinals cofinal in $\kappa=\aleph_\omega$. Also suppose $\Lambda>\kappa$ and the bounding number of $\left(\prod_{i<\omega}\lambda_i,<^*\right)$ is at least $\Lambda$. If $D$ is an ultrafilter on a set $X$ and $D$ is $(\theta,\kappa)$-regular for some $\theta<\kappa$, then $\chi_{<\kappa}(D)\ge\Lambda$.
\end{lemma}

We only need this lemma for the case $X=\kappa$, but we feel it better to distinguish between the underlying set $X$ and the regularity parameter $\kappa$.

Granted these two lemmas, we can prove $\uU^*(\aleph_\omega)=\aleph_{\omega+2}$ in Gitik's model as follows. Let $D$ be a uniform ultrafilter on $\aleph_\omega$. By Lemma~\ref{lem:dichotomy} it is either $(\omega,\kappa)$-regular or Rudin--Keisler above a strongly uniform ultrafilter $U$. In the latter case it is easy to see that $\chi^*(D)\geq\chi^*(U)$, and $\chi^*(U)=\aleph_{\omega+2}$ by Gitik's result, so $\chi^*(D)=\aleph_{\omega+2}$. In the former case, by standard facts from pcf theory there is some increasing sequence $(\lambda_i:i<\omega)$ of regular cardinals below $\aleph_\omega$ that carries a scale of length $\aleph_{\omega+2}$ (Gitik \cite{gitik2019spectrum} also gave such a set explicitly), so in particular the bounding number of $\left(\prod_{i<\omega}\lambda_i,<^*\right)$ is $\aleph_{\omega+2}$, and we apply Lemma~\ref{lem:regular-bound} with $\Lambda=\kappa^{++}$, $X=\kappa$ and $\theta=\omega$.

It remains to prove the two lemmas.

\begin{proof}[Proof of Lemma~\ref{lem:dichotomy}]
Our proof resembles that of Kanamori \cite[Theorem 2.5]{kanamori1976weakly}. First note that whenever a uniform ultrafilter $W$ on $\kappa$ is not strongly uniform, we can project it as follows. The sequence $(\aleph_i:i<\omega)$ fails to witness strong uniformity, so there exists $A\in W$ for which $|A\cap\aleph_i|<\aleph_i$ for all large enough $i<\omega$; by discarding an initial segment of $A$, we may assume this is true for every $i$, and moreover $A\cap\aleph_0=\emptyset$. Take $\aleph_{-1}=0$. We can define an injection $a:A\rightarrow\kappa$ that maps $A\cap[\aleph_i,\aleph_{i+1})$ into $[\aleph_{i-1},\aleph_i)$ for each $i$; note that $a(x)<x$. Extend $a$ to $e:\kappa\rightarrow\kappa$ by letting $e$ be identity outside $A$; it has the same projection effect as $a$.

Now suppose $D$ is uniform and does not project to any strongly uniform ultrafilter on $\kappa$. Recursively for $\xi<\kappa$ define $f_\xi:\kappa\rightarrow\kappa$ and $W_\xi$ as follows:
\begin{enumerate}
\item $f_0=\mathrm{id}$.
\item $W_\xi=(f_\xi)_*D$.
\item Apply the above construction to $W_\xi$ to get $A_\xi\in W_\xi$, $a_\xi:A_\xi\rightarrow\kappa$ and $e_\xi:\kappa\rightarrow\kappa$, and let $f_{\xi+1}=e_\xi\circ f_\xi$.
\item For limit $\eta$, let $f_\eta(x)$ be the eventual value of the non-increasing sequence $(f_\xi(x):\xi<\eta)$.
\end{enumerate}

We need to verify that $W_\xi$ is uniform in order for the recursion to go through. For $\xi<\kappa$ let $H_\xi(x)=\{\eta<\xi:f_\eta(x)\neq f_{\eta+1}(x)\}$ be the finite set of stages where the value at $x$ changes. We claim that $f_\xi(x)$ and $H_\xi(x)$ together recover $x$; this is because the maps $a_\eta$ are all injective. Thus uniformly for $y<\kappa$ we have
\[
|f^{-1}_\xi(y)|\leq|[\xi]^{<\omega}|\leq\max\{\aleph_0,|\xi|\}<\kappa.
\]
Consequently, if $C\subseteq\kappa$ has size strictly smaller than $\kappa$ then so does $f^{-1}_\xi[C]$. This shows $W_\xi=(f_\xi)_*D$ is uniform.

Now for $\xi<\kappa$ let $X_\xi=f^{-1}_\xi[A_\xi]\in D$. Note that $x\in X_\xi$ if and only if $f_\xi(x)>f_{\xi+1}(x)$, so no $x$ can belong to infinitely many $X_\xi$'s. The sequence $(X_\xi:\xi<\kappa)$ then witnesses $D$ is $(\omega,\kappa)$-regular.
\end{proof}

\begin{proof}[Proof of Lemma~\ref{lem:regular-bound}]
Our proof resembles Solomon's proof of $\mathfrak{b}\leq\uU$ \cite{solomon1977families}. Suppose $\Bcal\subseteq D$ has size smaller than $\Lambda$. We are going to find a set $Y\subseteq X$ that splits every $B\in \Bcal$, i.e., $|B\cap Y|\geq\kappa$ and $|B\cap Y^c|\geq\kappa$, so $\Bcal$ cannot be a base modulo $[X]^{<\kappa}$.

By assumption there is a sequence $(X_\alpha:\alpha<\kappa)$ such that each $X_\alpha$ is in $D$ and the intersection of any $\theta$ of them is empty; in other words, for each $x\in X$ the set $r(x):=\{\alpha<\kappa:x\in X_\alpha\}$ has size strictly smaller than $\theta$. Note that every set $A\in D$ has size at least $\kappa$.

Temporarily fix $B\in\Bcal$ and $n<\omega$. The set $B\cap X_\alpha$ has size at least $\kappa=\aleph_\omega$, so we can choose $E_\alpha\subseteq B\cap X_\alpha$ with size $\aleph_n$. Let $Z_\alpha=\bigcup_{x\in E_\alpha}r(x)$, so $|Z_\alpha|\leq \aleph_n\cdot\theta<\aleph_\omega$. Let $(\lambda_i:i<\omega)$ be as in the assumption; we may assume $\lambda_0>\theta$. For $i>n$, we can choose $h(i)<\lambda_i$ above every ordinal in $\bigcup_{\alpha<\lambda_{i-1}} Z_\alpha\cap\lambda_i$.

As $B\in\Bcal$ and $n<\omega$ vary, we get fewer than $\Lambda$ many such functions $h$, so by assumption they are $<^*$-bounded by some $g\in\prod_{i<\omega}\lambda_i$; we may assume $\lambda_{i-1}<g(i)<\lambda_i$ for $i>0$, so in particular $g$ is strictly increasing. Define
\[
Y=\{x\in X:g^{-1}[r(x)]\text{ has a largest element and it is even}\}.
\]

Let $E^{B,n}_\alpha$, $Z^{B,n}_\alpha$ and $h^{B,n}$ be as above. Take $N>n$ such that $h^{B,n}(i)<g(i)$ for all $i\geq N$. We claim that for all $i\geq N$ and $x\in E^{B,n}_{g(i)}$ we have $\max g^{-1}[r(x)]=i$. Clearly $i\in g^{-1}[r(x)]$ because $E^{B,n}_\alpha\subseteq X_\alpha$. Now suppose $j>i$ and $g(j)\in r(x)$; we have $g(i)<\lambda_i\leq\lambda_{j-1}$, so $g(j)\in Z^{B,n}_{g(i)}\cap\lambda_j$ and $g(j)<h^{B,n}(j)$, a contradiction.

Considering the even and odd $i>N$ respectively and $x\in E^{B,n}_{g(i)}$, we see that both $B\cap Y$ and $B\cap Y^c$ have size at least $\aleph_n$ for any $n$, hence at least $\aleph_\omega$, as desired.
\end{proof}

\noindent\textbf{AI disclosure:} no human intelligence was involved in writing this paper. All ideas are due to GPT 6 Astra.

\bibliographystyle{plain}
\bibliography{ref}

\end{document}